\documentclass[12pt]{article}
\usepackage[english,frenchb]{babel}
\usepackage[latin1]{inputenc}
\usepackage{amsfonts}
\usepackage{amsmath}
\usepackage{amsthm}
\usepackage{fullpage}

\begin{document}

\newcommand{\ABS}[1]{{\left| #1 \right|}} 
\newcommand{\PAR}[1]{{\left(#1\right)}} 
\newcommand{\SBRA}[1]{{\left[#1\right]}} 
\newcommand{\BRA}[1]{{\left\{#1\right\}}} 
\newcommand{\NRM}[1]{{\Vert #1\Vert}} 
\newcommand{\ind}{\mathrm{1}\hskip -3.2pt \mathrm{I}}

\newcommand{\moyn}[1]{{\left<#1\right>^{(n)}}} 

\newcommand{\gb}[1]{{\mathbb{E}\PAR{#1 e^{\beta \sum_{k=1}^n g(k,S_k)}}}} 
\newcommand{\gbi}[1]{{\mathbb{E}\PAR{#1 e^{\beta \sum_{i=1}^n  g(i,S_i)}}}}
\newcommand{\gbtilde}[1]{{\mathbb{E}\PAR{#1 e^{\beta \sum_{k=1}^n \tilde g(k,S_k)}}}} 
\newcommand{\gbtildei}[1]{{\mathbb{E}\PAR{#1 e^{\beta \sum_{i=1}^n \tilde g(i,S_i)}}}} 

\newcommand{\expgb}{ e^{\beta \sum_{k=1}^n g(k,S_k) } }

\newtheorem{ethm}{Theorem}[section]
\newtheorem{fthm}{Théorème}[section]
\newtheorem{ecor}[ethm]{Corollary}
\newtheorem{fcor}[fthm]{Corollaire}
\newtheorem{eprop}[ethm]{Proposition}
\newtheorem{fprop}[fthm]{Proposition}
\newtheorem{elem}[ethm]{Lemma}
\newtheorem{flem}[fthm]{Lemme}
\newtheorem{edefi}[ethm]{Definition}
\newtheorem{fdefi}[fthm]{Définition}
\newtheorem{erem}[ethm]{Remark}
\newtheorem{frem}[fthm]{Remarque}
\newtheorem{fex}{Exemple}
\newtheorem{eex}{Example}
\newtheorem{fhyp}[fthm]{Hypothèse}

\newcommand{\al}{\alpha}
\newcommand{\be}{\beta}
\newcommand{\de}{\delta}
\newcommand{\De}{\Delta}
\newcommand{\ga}{\gamma}
\newcommand{\Ga}{\Gamma}
\newcommand{\ep}{\epsilon}
\newcommand{\va}{\varphi}
\newcommand{\ka}{\kappa}
\newcommand{\la}{\lambda}
\newcommand{\La}{\Lambda}
\newcommand{\te}{\theta}
\newcommand{\Te}{\Theta}
\newcommand{\om}{\omega}
\newcommand{\Om}{\Omega}
\newcommand{\si}{\sigma} 
\newcommand{\ot}{\otimes}
\newcommand{\ti}{\times}
\newcommand{\cd}{CD(\rho ,\infty)}
\newcommand{\varep}{\varepsilon}
\newcommand{\na}{\nabla}
\newcommand{\ph}{\Phi}

\newcommand{\dA}{\mathbb{A}}
\newcommand{\dB}{\mathbb{B}}
\newcommand{\dC}{\mathbb{C}}
\newcommand{\dD}{\mathbb{D}}
\newcommand{\dE}{\mathbb{E}}
\newcommand{\dF}{\mathbb{F}}
\newcommand{\dG}{\mathbb{G}}
\newcommand{\dH}{\mathbb{H}}
\newcommand{\dI}{\mathbb{I}}
\newcommand{\dJ}{\mathbb{J}}
\newcommand{\dK}{\mathbb{K}}
\newcommand{\dL}{\mathbb{L}}
\newcommand{\dM}{\mathbb{M}}
\newcommand{\dN}{\mathbb{N}}
\newcommand{\dO}{\mathbb{O}}
\newcommand{\dP}{\mathbb{P}}
\newcommand{\dQ}{\mathbb{Q}}
\newcommand{\dR}{\mathbb{R}}
\newcommand{\dS}{\mathbb{S}}
\newcommand{\dT}{\mathbb{T}}
\newcommand{\dU}{\mathbb{U}}
\newcommand{\dV}{\mathbb{V}}
\newcommand{\dW}{\mathbb{W}}
\newcommand{\dX}{\mathbb{X}}
\newcommand{\dY}{\mathbb{Y}}
\newcommand{\dZ}{\mathbb{Z}}

\newcommand{\cA}{{\mathcal A }}
\newcommand{\cB}{{\mathcal B }}
\newcommand{\cC}{{\mathcal C }}
\newcommand{\cD}{{\mathcal D }}
\newcommand{\cE}{{\mathcal E }}
\newcommand{\cF}{{\mathcal F }}
\newcommand{\cG}{{\mathcal G }}
\newcommand{\cH}{{\mathcal H }}
\newcommand{\cI}{{\mathcal I }}
\newcommand{\cJ}{{\mathcal J }}
\newcommand{\cK}{{\mathcal K }}
\newcommand{\cL}{{\mathcal L }}
\newcommand{\cM}{{\mathcal M }}
\newcommand{\cN}{{\mathcal N }}
\newcommand{\cO}{{\mathcal O }}
\newcommand{\cP}{{\mathcal P }}
\newcommand{\cZ}{{\mathcal Z }}

\newenvironment{ithm}
 {\noindent \\ \noindent{\emph{\textbf{Théorème.}}} \begin{it}}
 {\end{it}}            
\newenvironment{iprop}
 {\noindent {\emph{\textbf{Proposition.}}} \begin{it}}
 {\end{it}\\}     


\newenvironment{eproof}
               {\noindent {\textbf{Proof.}}}

\newenvironment{fproof}
               {\noindent {\emph{\textbf{Preuve.}}}}

\newenvironment{enumeratei}{
 \renewcommand{\theenumi}{\roman{enumi}}
 \renewcommand{\labelenumi}{(\theenumi)} 
 \begin{enumerate}}{\end{enumerate}}


\title{Upper bound of a volume exponent for directed polymers in a random
 environment }

\author{Olivier \sc{Mejane}
        \\ \emph{Laboratoire de Statistiques et Probabilités,}
        \\ \emph{Université Paul Sabatier,
118 route de Narbonne}
     \\ \emph{31062 Toulouse France}   
     \\ email: Olivier.Mejane@lsp.ups-tlse.fr}

\maketitle

\section{Introduction}

The model of directed polymers in a random environment was introduced by
Imbrie and Spencer \cite{imbrie-spencer}. We focus here on a particular
model studied by Petermann \cite{markus} in his thesis~:
\\
let $(S_n)_{n \geq 0}$ be a random walk in $\dR ^d $ starting from the origin,
 with independent 
$\cN(0,I_d)$-increments, defined on a probability space $(\Om, \cF, \dP)$,
and let $g=(g(k,x), k \geq 1, x \in \dR ^d)$ be a stationary centred Gaussian
process with covariance matrix
$$cov(g(i,x),g(j,y))=\de _{ij} \Ga(x-y) \, ,$$
where $\Ga$ is a bounded integrable function on $\dR ^d .$
We suppose that this random media $g$ is defined on a probability space $(\Om
^g, \cG, P)$, where $(\cG_n)_{n \geq 0}$ is the natural filtration: 
$$\cG_n=\si \PAR{ g(k,x), 1 \leq k \leq n, x \in \dR ^d }$$ for $n \geq 1$ ($\cG_0$ 
being the trivial $\si$-algebra).

We define the Gibbs measure $\left< .\right>^{(n)}$ by~:
$$\moyn{f}=\frac{1}{Z_n} \gb{ f(S_1,\ldots,S_n) }$$
for any bounded function $f$ on $\PAR{\dR^d}^n$,
where $\be >0$ is a fixed parameter and $Z_n$ is the partition function~:
$$Z_n= \dE \PAR {\expgb}  \, .$$

Following Piza \cite{piza} we define the volume exponent
$$ \xi=\inf \BRA {\al >0 : \moyn{ \ind_{\BRA{\max_{k\leq n}\ABS{S_k} \leq
        n^{\al}}} } \xrightarrow[n \to \infty]{} 1 \quad \mbox{in} \quad P- \mbox{probability}} \, .$$
Here and in the sequel, $\ABS{x}=\max_{1 \leq i \leq d} \ABS{x_i}$ for any
$x=(x_1, \ldots , x_d) \in \dR^d$.
\\
Petermann obtained a result of superdiffusivity in dimension one,
in the particular case where $\Ga(x)=\frac{1}{2\la} e^{-\la \ABS{x}}$ for some
$\la>0$:
he proved that $ \xi \geq \frac{3}{5} $ for all $\be >0$
(for another result of superdiffusivity, see \cite{licea-newman-piza}).

Our main result gives on the contrary an upper bound for the volume 
exponent, in all dimensions~: 
\begin{equation} \label{upperbound}
\forall d \geq 1 \, , \forall \be >0 \quad  \xi \leq \frac{3}{4} \, .
\end{equation}

This paper is oragnized as follows:
\begin{itemize}
\item
In Section 2, we first extend exponential inequalities concerning independent
 Gaussian variables, proved by Carmona and Hu \cite{carmona-hu}, to the case
 of a stationary Gaussian process.Then, following Comets, Shiga and Yoshida
\cite{comets}, we combine these inequalities with martingale
methods and obtain a concentration inequality.

\item
In section 3, we obtain an upper bound for $d=1$, when we consider only the
value of the walk $S$ at time $n$, and not the maximal
one before $n$. 
\item
Section 4 is devoted to the proof of (\ref{upperbound}). 

\end{itemize}

\section{Preliminary: A concentration inequality}
 \subsection{Exponential inequalities}
\begin{elem} \label{expo_ineq}
Let $(g(x),x \in \dR^d)$ be a family of Gaussian centered random variables with
common variance $\si^2 > 0$.
We fix $q,\be >0$,$(x_1, \ldots , x_n) \in \PAR{\dR^d}^n$ and $(\la_1, \ldots , \la_n)$ in $\dR^n$.
Then for any probability measure $\mu$ on $\dR^d$~:

$$e^{-\frac{\be^2 \si^2}{2} q } \leq \dE \PAR{ \frac{ e^{\be \sum_{i=1}^n
      \la_i g(x_i)}}{\PAR{\int_{\dR} e^{\be g(x) } \mu(dx)}^q}} \leq 
  e^{\frac{\be^2 \si^2}{2} \PAR{q+\sum_{i=1}^n \ABS{\la_i}}^2}$$
\end{elem}

The proof is identical to the one made by Carmona and Hu in a discrete
framework ($\mu$ is the sum of Dirac masses),and is therefore 
omitted.

\begin{elem} \label{Carm_Hu}
Let $(g(x),x \in \dR^d)$ be a centred Gaussian process with covariance matrix

$cov(g(x),g(y))=\Ga(x-y)$. Let $\si^2=\Ga(0)$, and let $\mu$ be a probability
measure on $\dR^d$.
Then for all $\be >0$, there are non-negative constants $c_1=c_1(\be,\si^2)$ and $c_2=c_2(\be,\si^2)$
such that~:
$$
-c_1 \int\!\!\!\int_{\dR^d}\Ga(x-y) \,\mu(dx) \,\mu(dy)
\leq E \PAR {\log \int_{\dR} e^{\be g(x)-\frac{\be^2 \si^2}{2}} \, \mu(dx)}
\leq -c_2 \int\!\!\!\int_{\dR^d}\Ga(x-y) \,\mu(dx) \,\mu(dy) \, .
$$
In particular,
$$ -c_1 \si^2 \leq   
 E \PAR {\log \int_{\dR^d} e^{\be g(x)-\frac{\be^2 \si^2}{2}}
 \, \mu(dx)}
\leq 0  \, .
$$
\end{elem}

\begin{eproof}
 Let $\{B_x(t), t \geq 0\}_{x \in \dR^d}$ be the family of centred Gaussian processes such
 that $$E \PAR { B_x(t) B_x(s) }= \inf(s,t) \, \Ga(x-y) \, ,$$ with $B_x(0)=0$
for all $x \in \dR^d$.
 Define
 $$ X(t)=\int_{\dR^d}M_x(t) \mu(dx) , \quad
   t \geq 0,
 $$
where $M_x(t)=e^{\be B_x(t)-\frac{\be^2 \si^2 t}{2}}$.
Since $dM_x(t)=\be M_x(t) dB_x(t)$, one has
$$d\left<M_x,M_y \right>_t=\be^2 M_x(t) M_y(t) d\left<B_x,B_y\right>_t=
\be^2 e^{\be(B_x(t)+B_y(t))-\be^2 \si^2 t} \Ga(x-y) dt \, ,$$ and
$d\left<X,X\right>_t=\int\!\!\!\int_{\dR^{d}}\be^2 e^{\be(B_x(t)+B_y(t))-\be^2 \si^2 t} \, \Ga(x-y)
\mu(dx) \mu(dy)\,
 dt.$

Thus,
by Ito's formula,

$$ E \PAR {\log X_1}=-\frac{\be^2}{2}
\int\!\!\!\int_{\dR^2} \mu(dx) \,\mu(dy) \Ga(x-y) 
\int_0^1 E \PAR {\frac{e^{\be(B_x(t)+B_y(t))-\be^2 \si^2 t}}{X_t^2}} dt \, .$$

By Lemma \ref{expo_ineq}, we have for all $t$~:

$$ e^{-\be^2 \si^2 t}\leq 
E \PAR {\frac{e^{\be(B_x(t)+B_y(t))-\be^2 \si^2 t}}{X_t^2}}
=E \PAR {\frac{e^{\be(B_x(t)+B_y(t))}}
{\PAR{\int_{\dR} e^{\be B_x(t)} \mu(dx)}2}} \leq e^{8\be^2 \si^2 t} \, .$$

Hence~:

$$ -\frac{e^{8\be^2 \si^2} -1}{16\si^2} 
\int\!\!\!\int_{\dR^d} \Ga(x-y)  \mu(dx) \,\mu(dy) \leq
E \PAR {\log X_1}  \leq
-\frac{1-e^{-\be^2 \si^2}}{2\si^2} 
\int\!\!\!\int_{\dR^d} \Ga(x-y)  \mu(dx) \,\mu(dy) \, ,$$

which concludes the proof since $ X_1 \stackrel{d}{=}\int_{\dR^d} e^{\be
  g(x)-\frac{\be^2 \si^2}{2}} \, \mu(dx).  $ 

\end{eproof}

\subsection{A concentration result}

\begin{eprop} \label{concentration}
Let $\nu > \frac{1}{2}$.
For $n \in \dN$, $j \leq n$ and $f_n$ a bounded function, we note
$W_{n,j}=\gb{f_n(S_j)}$. Then for $n \geq n_0(\be,\si^2,\nu)$,

$$P \PAR { \ABS{ \log W_{n,j}
 - E \PAR{\log W_{n,j}} } \geq n^{\nu}} \leq \exp(-\frac{1}{4} n^{(2\nu-1)/3})
\, .$$
\end{eprop}

\begin{eproof}
Following Comets, Shiga and Yoshida (\cite{comets},Proposition 1.5),we use a large deviation estimate for sum
 of martingale-differences, which is a slight extension of a result of Lesigne and Voln\'y
(\cite{lesigne}, Theorem 3.2) :

\begin{elem}
Let $(X_n^i,1 \leq i \leq n)$ a sequence of martingale-differences and
let $S_n=\sum_{i=1}^n X_n^i$.
 Suppose that there exists $K>0$ such that $\dE (e^{\ABS{ X_n^i}}) \leq K$ for all
$i$ and $n$.
Then for any $\nu > \frac{1}{2}$, and for $n \geq n_0(K,\nu) $,

$$\dP \PAR{ |S_n|> n^{\nu} } \leq \exp(-\frac{1}{4} n^{(2\nu-1)/3}) \, . $$ 
\end{elem}

In our case we define 
$X_{n,j}^i=E \PAR{ \log W_{n,j} |\cG_i} -E \PAR{ \log W_{n,j} |\cG_{i-1}}$
so that 
$$\log(W_{n,j}) - E (\log W_{n,j})=\sum_{i=1}^n X_{n,j}^i \, .$$
It is sufficient to prove that there exists $K>0$ such that
$\dE(e^{\ABS{X_{n,j}^i}})
\leq K$
for all $i$ and $(n,j)$.

For this, we introduce: 

$$e_{n,j}^i= f_n(S_j) \exp \PAR{\sum_{1 \leq k \leq n, k\neq i} \be g(k,S_k)}, \, W_{n,j}^i=E(e_{n,j}^i).$$

If $E_i$ is the conditional expectation with respect to $\cG_i$,
 then $E_i \PAR{\log W_{n,j}^i}=E_{i-1}\PAR{\log W_{n,j}^i},$  so that:

\begin{equation} \label{0step}
X_{n,j}^i=E_i \PAR {\log Y_{n,j}^i } 
-E_{i-1} \PAR {\log Y_{n,j}^i } \, ,
\end{equation}
with 
 \begin{equation} \label{rewriting}
Y_{n,j}^i =e^{-\be^2/2}\frac{W_{n,j}}{W_{n,j}^i}= \int_{\dR^d} e^{\be g(i,x)-\be ^2/2} \mu_{n,j}^i(dx) \, ,
\end{equation}

$\mu_{n,j}^i$ being the random probability measure:

$$\mu_{n,j}^i(dx)=\frac{1}{W_{n,j}^i} \dE \PAR{e_{n,j}^i | S_i=x } \dP (S_i\in
dx) \, .$$

Since $\mu_{n,j}^i$ is measurable with respect to 
$\cG_{n,i}= \si \PAR {g(k,x), 1 \leq k \leq n, k \neq i, x \in \dR^d)}$,

we deduce from Lemma \ref{Carm_Hu} that there exists a constant $c=c(\be,\si^2)>0$, which does not depend on $(n,j,i)$ , such that~:

$$-c \leq E \PAR{ \log \int_{\dR^d}  e^{\be g(i,x)-\be ^2/2} \mu_{n,j}^i(dx)
 \, | \, \cG_{n,i}} \leq 0 \, ,$$

and since $\cG_{i-1} \subset \cG_{n,i}$, we obtain:

\begin{equation} \label{1step}
 0 \leq -E_{i-1} \PAR {\log Y_{n,j}^i} \leq c \, ,
\end{equation}

Thus we deduce from (\ref{0step}) and (\ref{1step}) that for all 
$\te \in \dR$ :

$$E \SBRA{ e^{\te X_{n,j}^i}} \leq 
e^{c\te^+} E \SBRA { e^{ \te E_i \PAR { \log Y_{n,j}^i}}}$$

with $\te^+:=\max(\te,0)$.

By Jensen's inequality,

$$e^{ \te E_i \PAR { \log Y_{n,j}^i}}  \leq
 E_i \SBRA {\PAR{ Y_{n,j}^i}^{\te}}$$

so that

$$E \SBRA{ e^{\te X_{n,j}^i}} \leq
e^{c\te^+} E \SBRA{ \PAR {Y_{n,j}^i}^{\te}}= 
e^{c\te^+} E  \SBRA {E \SBRA {\PAR{ Y_{n,j}^i}^{\te}|\cG_{n,i}}} \, ,$$

Assume now that $\te \in \{-1,1\}$, hence in both cases, the function
 $x \, \mapsto \, x^{\te}$ 
is convex; using (\ref{rewriting}), we obtain
$\PAR{Y_{n,j}^i}^{\te} \leq \int_{\dR^d}
e^{\te(\be g(i,x)-\be ^2/2)} \mu_{n,j}^i(dx) \, ,$so that:

\begin{eqnarray*}
E \SBRA {\PAR{ Y_{n,j}^i}^{\te}| \cG_{n,i}} &\leq& 
\int_{\dR^d} E \PAR{ e^{\te(\be g(i,x)-\be ^2/2)}| \cG_{n,i}} \mu_{n,j}^i(dx) \\
 &=& \int_{\dR^d} E \PAR{ e^{\te(\be g(i,x)-\be ^2/2)}} \mu_{n,j}^i(dx) \\
 &=& E \PAR{ e^{\te(\be g(1,0)-\be ^2/2)}} \, ,
\end{eqnarray*}

using that $g(i,x)$ is independent from
$\cG_{n,i}$, and is distributed as $g(1,0)$
for all $i$ and $x$.

We conclude that for all $n$ and $1 \leq i , j \leq n$,

\begin{eqnarray*}
E \SBRA{ e^{\ABS{ X_{n,j}^i}}} & \leq &  
 E \SBRA{ e^{ X_{n,j}^i}} +  E \SBRA{ e^{- X_{n,j}^i}} \\
 &\leq & K := e^{c}+ e^{\be^2} \, .  
\end{eqnarray*}

\end{eproof}

\begin{ecor} \label{cor_concentration}
Let $\nu > \frac{1}{2}$. Then  $P$-almost surely, there exists $n_0$ such
that for every $n \geq n_0$, every $ j\leq n$ and every Borel set $B(j,n)$,

$$ \ABS{ \log \moyn {\ind_{S_j \in B(j,n)}} 
 -E \PAR{\log \moyn{\ind_{S_j \in B(j,n)}} }} \leq 2n^{\nu} $$
\end{ecor}

\begin{eproof}
Let us write $A_{n,j}=\BRA {\ABS{ \log \gb{f_n(S_j)}
 - E \SBRA{\log \gb{f_n(S_j)}} } \geq n^{\nu}}$.
Proposition \ref{concentration} implies that  
$$P \PAR{ \cup_{j\leq n} A_{n,j}} \leq n \exp(-\frac{1}{4} n^{(2\nu-1)/3}) \, .$$
Hence, by Borel-Cantelli, $P$-almost-surely there exists $n_0$ such
that for every $n \geq n_0$ and every $j \leq n$:
 $$
\ABS{ \log \gb{f_n(S_j)}
 - E \SBRA{\log \gb{f_n(S_j)} }} \leq n^{\nu} \, .
$$
Then one applies  this result to 
$f_n(x)=\ind_{x \in B(j,n)}$ and to $f_n(x)=1$.

\end{eproof}

\section{A first result}

\begin{ethm} \label{main}
If $d=1$, for all $\al >\frac{3}{4}$, 
$$ \moyn{ \ind_{\ABS{S_n} \geq n^{\al}} } \xrightarrow [n \to \infty]{P-a.s.} 0$$
\end{ethm}

\begin{eproof}
We first prove:
\begin{eprop} \label{mean_control}
For all $n \geq 0$,
$$E \PAR{\log \moyn {\ind_{ S_n \geq n^{\al} } }}
 \leq  -\frac{1}{2} n^{2\al -1} \, .$$
\end{eprop}
\begin{eproof}
Let us fix $a >0$ and $\la>0$.
We first write~:
$$\moyn{\ind_{S_n \geq a}} \leq e^{-\la a} \moyn{ e^{\la S_n}}
= e^{-\la a + n\la^2 /2} \moyn{ M_n^{\la}} \, ,$$
$(M_n^{\la}:=e^{\la S_n - n\la^2 /2})_{n \geq 0}$ being a nonnegative
 martingale. 
Then we have by Girsanov's Theorem~:
\begin{eqnarray}
 \moyn{ M_n^{\la}}&=&\frac{\gb{ M_n^{\la}}}{Z_n} \\
   &=&\frac{\dE \PAR { e^{\be \sum_{k=1}^n g(k,S_k+k\la)} }}{Z_n} \\
    &=&\frac{\dE \PAR{ e^{\be \sum_{k=1}^n g^{\la}(k,S_k) }}} 
                          {\dE \PAR{ e^{\be \sum_{k=1}^n g(k,S_k) }}} \, ,
\end{eqnarray}
where in the last equation we denote by $g^{\la}$ the translated
 environment 
$$(g^{\la}(k,x):=g(k,x+k\la), k \geq 1, x \in \dR^d).$$
By stationarity, this environment has the same distribution 
as $(g(k,x)), k \geq 1, x \in \dR^d)$, hence
$$ \dE{ e^{\beta \sum_{k=1}^n g^{\la}(k,S_k)}}\stackrel{d}{=}
 \dE{ e^{\beta \sum_{k=1}^n g(k,S_k)}} \, ,$$
and thus 
$$E \PAR {\log \moyn{ M_n^{\la}}}=0 \, .$$

We conclude that 
$$E \PAR {\log \moyn{1_{S_n \geq a}} }\leq -\la a + n\la^2 /2$$
and  by taking $\la=\frac{a}{n}$ we obtain the upper bound $-\frac{a^2}{2n}$,
which gives the result when $a=n^{\al}$.
\end{eproof}

Assume now that $\nu > \frac{1}{2}$ and $\al > \frac{\nu+1}{2}$.
We deduce from the last proposition  and Corollary \ref{cor_concentration} that
$P$-almost-surely for large $n$:
$$\log \moyn{\ind_{{S_n} \geq n^{\al}} }\leq 
 -\frac{1}{2} n^{2\al -1} + 2 n^{\nu} \xrightarrow[n \to + \infty]{} - \infty$$

Since this is true for all $\nu > \frac{1}{2}$, 
$$\moyn{ \ind_{S_n \geq n^{\al}}} \xrightarrow[n \to + \infty]{P-a.s.} 0$$
for all $\al > \frac{3}{4}$.

But we have in the same way that $$\moyn{ \ind_{S_n \leq -n^{\al}}}
\xrightarrow[n \to + \infty]{P-a.s.} 0$$
for all $\al > \frac{3}{4}$, which ends the proof of Theorem \ref{main}.

\end{eproof}

\section{Extension to the maximum}

We now extend the previous result to the maximal deviation 
 from the origin and to all dimensions $d$~:

\begin{ethm} \label{max}
For all $d \geq 1$ and $\al >\frac{3}{4}$,  
$$ \moyn{ \ind_{\BRA{\max_{k\leq n}\ABS{S_k} \geq n^{\al}}} } \xrightarrow [n \to
\infty]{P a.s.} 0 \, .$$
\end{ethm}

\begin{eproof}
\\
We will use the following notations:
for $x \in \dR^d$ and $r \geq 0$, $B(x,r)=\{ y\in \dR^d, \ABS{y-x}\leq r\}$.
For $\al \geq 0$ and $j=(j_1, \ldots, j_d) \in \dZ^d$,
$B_j^{\al}=B(jn^{\al},n^{\al})$.
We will use the fact that the union of the balls 
$(B_j^{\al}, j \in \PAR{2\dZ}^d \backslash \BRA{0})$ form a partition of 
$\dR^d  \backslash B(0,n^{\al})$.

We first need the same kind of upper bound as in Proposition
\ref{mean_control}:

\begin{eprop} 
Let $n \geq 0$ and $k \leq n$. Then for any $j\in \dZ^d$ and $\al >0$,
$$E \PAR{\log \moyn {\ind_{ S_k \in B_j^{\al}} } }
 \leq  \frac{-n^{2\al -1}}{2} \sum_{i=1}^d (j_i-\ep_i)^2 \, ,$$
where $\ep_i=sgn(j_i) (=0  \, \mbox{if} \, j_i=0)$.

\end{eprop}

\begin{eproof}

Let note $a_{k,j}^{\al}=\gbi{\ind_{ S_k \in B_j^{\al} }}$, so that
$\moyn {\ind_{ S_k \in B_j^{\al} } } = \frac{a_{k,j}^{\al}}{Z_n}$.

Let be $\la=\frac{\tilde{\la}}{k}$ with 
$\tilde{\la}_i=(j_i-\ep_i)n^{\al} , \, 1 \leq i \leq d$;
then let us define the martingale
$$M_p^{\la,k}=
\left\{
\begin{array}{cc}
e^{\la.S_p - p\NRM{\la}^2 /2}  & if  \, p \leq k \\
e^{\la.S_k - k\NRM{\la}^2 /2}  & if  \, p > k ,
\end{array}
\right.
$$
where $x.y$ denotes the usual scalar product in $\dR^d$ and $\NRM{x}$
the associated euclidean norm.
Under the probability $\dQ^{\la,k}$ defined by Girsanov's change associated
 to this martingale, $(S_p)_{p \geq 0}$ has the law
 of the following shifted random walk under $\dP$:

$$\tilde S_p=
S_p + \tilde{\la} \PAR{\frac{p}{k} \wedge 1} \, .  
$$

It follows that:
\begin{equation} \label{girsanov}
a_{k,j}^{\al}= \gbtildei{ e^{\frac{-1}{k}(\tilde{\la}. S_k + \NRM{\tilde{\la}}^2/2) }
 \ind_{S_k \in  B_j^{\al}-\tilde{\la} } } \, ,
\end{equation}
where $\tilde{g}(i,x)=g(i,x+ \tilde{\la} (\frac{i}{k} \wedge 1))$.

Now we notice that on the event 
$\BRA{ S_k\in B_j^{\al}-\tilde{\la}}$, one has $\tilde{\la}.S_k \geq 0$:
indeed if we write $S_k=(S_k^1, \ldots, S_k^d)$, then for any $1 \leq i \leq d$, 
$\ABS{S_k^i -j_i n^{\la} + \tilde{\la}_i} \leq n^{\al}$, hence~:
\begin{itemize}
\item
for $j_i \geq 1$, $\tilde{\la}_i=j_i-1 \geq 0$ and 
$ 0 \leq S_k^i \leq 2 n^{\al}$,
\item
for $j_i \leq -1$, $\tilde{\la}_i=j_i+1 \leq 0$ and 
$ -2 n^{\al} \leq S_k^i \leq 0$,
\item
for $j_i=0$, $\tilde{\la}_i=0$,
\end{itemize}
so that in all cases $\tilde{\la}_i S_k^i \geq 0$ and thus 
$\tilde{\la}.S_k \geq 0$.

Therefore on the event $\BRA{ S_k\in B_j^{\al}-\tilde{\la}}$,
 $$e^{\frac{-1}{k}(\tilde{\la}. S_k + \NRM{\tilde{\la}}^2/2)} \leq
e^{\frac{-\NRM{\tilde{\la}}^2}{2n}}
\leq
 e^{\frac{-n^{2\al -1}}{2}\sum_1^d (j_i-\ep_i)^2 },$$
and (\ref{girsanov}) leads to~:
$$a_{k,j}^{\al} \leq e^{\frac{-n^{2\al -1}}{2}\sum_1^d (j_i-\ep_i)^2 }\gbtildei { \ind_{S_k\in B_j^{\al}-\tilde{\la}}}.$$

On the other hand, $Z_n \geq \gbi { \ind_{S_k \in B_j^{\al}-\tilde{\la}} }$,
and since by stationarity the environment $\tilde g$ has the same
distribution as $g$, it follows that for all $j \in \dZ^d$,

$$\dE \PAR{ \log \moyn {\ind_{ S_k \in B_j^{\al}} } } \leq  
\frac{-n^{2\al-1}}{2} \sum_{i=1}^d (j_i-\ep_i)^2  \, .$$

\end{eproof}

Let $\nu > \frac{1}{2}$.
We deduce from the last proposition and from Corollary \ref{cor_concentration}
(with $B(j,n)=B_j^{\al}$) that for $n \geq n_0$ and all 
$j \in \dZ^d$:

$$\log \moyn {\ind_{ S_k \in B_j^{\al}} } \leq 2 n^{\nu}- \frac{n^{2\al-1}}{2} \sum_{i=1}^d (j_i-\ep_i)^2  \, .$$

So, for $n \geq n_0$,
$\moyn{\ind_{ \ABS{S_k} \geq n^{\al} } } \leq \sum_{j \in (2\dZ)^d\backslash
  \BRA{0}} \, e^{2 n^{\nu}- \frac{n^{2\al-1}}{2} \sum_{i=1}^d (j_i-\ep_i)^2},$
 and~: 

$$\moyn{\ind_{\BRA{\max_{ k\leq n} \ABS{S_k} \geq  n^{\al}} } } 
\leq \sum_{k=1}^{n} \moyn{\ind_{ \ABS{S_k} \geq n^{\al} }}     
\leq \sum_{j \in (2\dZ)^d\backslash \BRA{0}} n  e^{2 n^{\nu}-
 \frac{n^{2\al-1}}{2} \sum_{i=1}^d (j_i-\ep_i)^2 } \, .$$

but by symmetry, for any $C >0$,
 $$\sum_{j \in (2\dZ)^d\backslash \BRA{0}  } e^{-C\sum_{i=1}^d (j_i-\ep_i)^2 }
\leq 2d  \sum_{ j_1 \geq 2} e^{-C(j_1-1)^2} \prod_{i=2}^d \sum_{ j_i \in 2\dZ} e^{-C(j_i-\ep_i)^2}$$

 and using that $\sum_{ j \geq 2} e^{-C(j-1)^2} \leq
 \sum_{j \geq 1} e^{-Cj}= \frac{ e^{-C}}{1-e^{-C}},$

we conclude that for some constant $C(d)>0$, and for $n \geq n_0$:

$$\moyn{\ind_{\BRA{\max_{ k\leq n} \ABS{S_k} \geq  n^{\al}} } } 
\leq
C(d) \,  n e^{2 n^{\nu}} \frac{ e^{-n^{2\al-1}/2}}{1-e^{-n^{2\al-1}/2}}.$$ 

Thus for all $\al > \frac{\nu+1}{2}$, $P$-a.s. ,

$$\moyn{\ind_{\BRA{\max_{ k\leq n} \ABS{S_k} \geq  n^{\al}} } } 
 \xrightarrow[n \to \infty]{} 0 \quad ,$$

and we conclude as in the proof of Theorem \ref{main}.

\end{eproof}
\nocite{sherrington,song-zhou,albeverio,bolthausen}
\bibliographystyle{plain}
\bibliography{biblio}

\begin{thebibliography}{10}

\bibitem{albeverio}
Sergio Albeverio and Xian~Yin Zhou.
\newblock A martingale approach to directed polymers in a random environment.
\newblock {\em J. Theoret. Probab.}, 9(1):171--189, 1996.

\bibitem{bolthausen}
Erwin Bolthausen.
\newblock A note on the diffusion of directed polymers in a random environment.
\newblock {\em Comm. Math. Phys.}, 123(4):529--534, 1989.

\bibitem{carmona-hu}
P.~Carmona and Y.~Hu.
\newblock On the partition function of a directed polymer in a gaussian random
  environment.
\newblock {\em Preprint}, 2001.

\bibitem{sherrington}
F.~Comets and J.~Neveu.
\newblock The {S}herrington-{K}irkpatrick model of spin glasses and stochastic
  calculus: the high temperature case.
\newblock {\em Comm. Math. Phys.}, 166(3):549--564, 1995.

\bibitem{comets}
F.~Comets, T.~Shiga, and N.~Yoshida.
\newblock Directed polymers in random environment: path localization and strong
  disorder.
\newblock {\em Preprint}, 2001.

\bibitem{imbrie-spencer}
J.~Z. Imbrie and T.~Spencer.
\newblock Diffusion of directed polymers in a random environment.
\newblock {\em J. Statist. Phys.}, 52(3-4):609--626, 1988.

\bibitem{lesigne}
Emmanuel Lesigne and Dalibor Voln{\'y}.
\newblock Large deviations for martingales.
\newblock {\em Stochastic Process. Appl.}, 96(1):143--159, 2001.

\bibitem{licea-newman-piza}
C.~Licea, C.~M. Newman, and M.~S.~T. Piza.
\newblock Superdiffusivity in first-passage percolation.
\newblock {\em Probab. Theory Related Fields}, 106(4):559--591, 1996.

\bibitem{markus}
M.~Petermann.
\newblock Superdiffusivity of directed polymers in random environment.
\newblock {\em Part of thesis}, 2000.

\bibitem{piza}
M.S.T. Piza.
\newblock Directed polymers in a random environment: Some results on
  fluctuations.
\newblock {\em Journal of Statistical Physics}, 89(3-4):581--603, 1997.

\bibitem{song-zhou}
Renming Song and Xian~Yin Zhou.
\newblock A remark on diffusion of directed polymers in random environments.
\newblock {\em J. Statist. Phys.}, 85(1-2):277--289, 1996.

\end{thebibliography}

\end{document}